\documentclass[preprint,12pt]{elsarticle}

\usepackage{amssymb}
\usepackage{bm}
\usepackage{amsfonts,amsmath,amsthm,amsxtra,amscd,amssymb}
\usepackage{verbatim,graphicx} 
\usepackage{tikz}

\def\a{\alpha}
\def\b{\beta}

\def\g{\gamma}
\def\G{\Gamma}

\def\L{\Lambda}

\def\o{\omega}

\def\s{\sigma}
\def\S{\Sigma}

\def\z{\zeta}

\def\P #1{\partial_{#1}}
\def\ch #1{{\rm Char}\ #1}
\def\bsy{\boldsymbol}
\def\wh{\widehat}

\newtheorem{thm}{Theorem}[section]
\newtheorem{prop}[thm]{Proposition}
\newtheorem{lem}[thm]{Lemma}
\newtheorem{cor}[thm]{Corollary}

\newtheorem{defn}{Definition}[section]

\newtheorem{rem}{Remark}[section]

\newproof{pf}{Proof}

\journal{Differential Geometry and its Applications}

\begin{document}
\baselineskip=0.17in
\begin{frontmatter}



\title{A Solvable String on a Lorentzian Surface}


\author[L1]{Jeanne N Clelland\fnref{J1}}\ead{Jeanne.Clelland@colorado.edu}
\fntext[J1]{The first author was supported in part by NSF grant DMS-1206272.}
\author[L2]{Peter J Vassiliou}\ead{peter.vassiliou@canberra.edu.au}

\address[L1]{Department of Mathematics, 395 UCB, University of Colorado Boulder, CO 80309-0395}
\address[L2]{Program in Mathematics \& Statistics, University of Canberra, ACT, Australia, 2601}


\begin{abstract}
It is shown that there are nonlinear sigma models which are Darboux integrable and possess a solvable Vessiot group in addition to those whose Vessiot groups are central extensions of semi-simple Lie groups. They govern harmonic maps between Minkowski space $\mathbb{R}^{1,1}$ and  certain complete, non-constant curvature 2-metrics. The solvability of the Vessiot group permits a reduction of the general Cauchy problem to quadrature. We treat the specific case of harmonic maps from Minkowski space into a non-constant curvature Lorentzian 2-metric, $\bsy{\lambda}$. Despite the completeness of $\bsy{\lambda}$ we exhibit a Cauchy problem with real analytic initial data which blows up in finite time. We also derive a hyperbolic Weierstrass representation formula for all harmonic maps from $\mathbb{R}^{1,1}$ into $\bsy{\lambda}$. 

\end{abstract}


\end{frontmatter}

\vskip 10 pt
\begin{center}
     \end{center}
\section{Introduction}
\label{Intro}
The linear wave equation
\begin{equation}\label{linwave}
\frac{\partial^2 u}{\partial t^2}=c^2\sum_{i=1}^n\frac{\partial^2 u}{\partial {x_i}^2}
\end{equation}
on $n$-dimensional Euclidean space is one of the most fundamental in all of mathematics and physics, modeling a host of phenomena, often to a first approximation, but sometimes exactly. At the simplest level, $n=1$, equation (\ref{linwave}) models the motion of a frictionless elastic string moving in the plane  under its elastic forces. In this case $u(t,x_1)$ represents the transverse displacement from a reference $x$-axis, fixed in the plane, of the point on the string located at rest position $x_1$. It can be shown that, contrary to what is found in many textbooks on this elementary topic, equation (\ref{linwave}) is {\it not} an approximate model  for small oscillations of the string but rather an {\it exact} physical model \cite{BleeckerCsordas}, \cite{ClellandVassiliou}.

It is reasonable to wonder about the problem of an elastic string constrained to vibrate on any Riemannian surface, such as, say, the surface of a sphere, or more abstractly, on a 2-manifold with a smooth Riemannian metric, such as the Poincar\'e half-plane.  To this end, let $M$ be a Riemannian surface equipped with smooth metric $g$, and let $\G:\mathbb{R}^2_{(\xi,\tau)}\to M$ be a smooth map such that for each $\tau$, $\G(\cdot, \tau):\mathbb{R} \to M$ is an immersion into $M$ whose image is a curve that models an elastic string vibrating without friction on $M$. Here $\tau$ represents time, and $\xi$ is a parameter along the curve (string). Let $\nabla$ denote the Levi-Civita connection for $(M,g)$. It can then be shown that the equation of motion of the string is given by \cite{ClellandVassiliou}
\begin{equation}\label{string}
\nabla_{\G_\tau}\G_\tau=c^2\nabla_{\G_\xi}\G_\xi,
\end{equation}
where subscripts such as the `$\tau$' in $\G_\tau$ denote partial differentiation with respect to $\tau$, and $c$ is a constant related to the tension in the string and its line density. Equation (\ref{string}) is the Euler-Lagrange equation corresponding to Lagrangian density
$$
{L}(\G)=
\frac{1}{2}g_{\a\b}(\G)\left(-\frac{1}{c^2}\frac{\partial \G^\a}{\partial \tau}\frac{\partial \G^\b}{\partial \tau}
+\frac{\partial \G^\a}{\partial \xi}\frac{\partial\G^\b}{\partial\xi}\right).
$$
We see that $\Gamma$ can be viewed as a {\it harmonic map} $\G:\mathbb{R}^{1,1}\to (M,g)$, where $\mathbb{R}^{1,1}$ is Minkowski space with metric
$$
\eta=d\xi^2-c^2\,d\tau^2.
$$
In this case $\mathcal{E}(\G)$ is the {\it energy} of $\G$:
\begin{equation}\label{energy}
\mathcal{E}(\G)=\frac{1}{2}\int\eta^{ij}g_{\a\b}\frac{\partial \G^\a}{\partial x_i}\frac{\partial \G^\b}{\partial x_j}\,d{\rm vol}_M,
\end{equation}
where we set $x_1=\xi, x_2=\tau$. 

Harmonic maps where the source manifold is a Minkowski space $\mathbb{M}$ with its usual flat Lorentzian metric are generally called {\it wave maps}. (In the 2-dimensional case $\mathbb{M}=\mathbb{R}^{1,1}$, they are sometimes called 1+1-wave maps.) There is a well-known geometric literature on wave maps that has developed over the last two decades, especially regarding their existence as solutions of completely integrable systems; see \cite{Guest} for a textbook account with many references. There is also a closely related physics literature where the relevant systems are known as {\em nonlinear sigma models}; see  \cite{Zak}. 

While (\ref{linwave}) is a linear partial differential equation and the equation of a string in the Euclidean plane is a linear problem, this is not the case for strings vibrating on a curved Riemannian---or, for that matter, pseudo-Riemannian---surface, despite the innocuous appearance of equation (\ref{string}). Indeed, if there are no elastic forces acting on the string, then the left-hand side of (\ref{string}) vanishes and the equation is that of a geodesic in $M$; this reflects the expected result that the equilibrium configuration of the string lies along a geodesic. As is well known, geodesic equations are nonlinear ordinary differential equations, and in general they cannot be exactly solved in terms of known functions. Thus, equation (\ref{string}) precisely generalises the geodesic equation of a Riemannian or pseudo-Riemannian manifold, and at the same time, it generalises the wave equation (\ref{linwave}). It derives its significance, in part, from this fact.

Thus, just as in the case of geodesic flow, a question naturally arises as to those Riemannian or pseudo-Riemannian surfaces for which the corresponding wave map equations governing the vibration of frictionless elastic strings exhibit ``integrable behaviour." Of course, this question is not new: it is well known that for a wide class of manifolds $M$, the corresponding equations are integrable in the sense of soliton theory, \cite{Guest}, \cite{TerngUhlenbeck}. But here we are interested in a stronger and more explicit form of integrability known as {\em Darboux integrability} \cite{IveyLandsbergBook}, \cite{Vassiliou2000}, \cite{AFV}, \cite{Zhiber}, \cite{ZhiberSokolov}. Differential systems satisfying this form of integrability often admit Weierstrass-type representations, and, as has recently been shown, their Cauchy problems can be reduced to ordinary differential equations of Lie type \cite{AndersonFels}, \cite{Vassiliou2012}, \cite{Bryant}, \cite{Carinena}, \cite{DoubrovKomrakov}.  In this context we note an ommision in the recent classification of Ream \cite{Ream} which aims to give a complete list of all wave maps systems into Riemannian surfaces which are Darboux integrable at order 2 or 3. A small error in a computation caused the author to miss an important case which we describe herein and which, we believe, completes the classification.  

Another intriguing question in relation to harmonic maps from split signature source manifolds is that of global existence for the Cauchy problem. In this case, Gu \cite{Gu} established global existence when the target manifold has a complete Riemannian metric. However, for split signature target metrics the global existence question appears to be more delicate. In that case Gu's theorem doesn't hold, as observed, for instance, in work of Terng and Uhlenbeck \cite{TerngUhlenbeck}. For this reason we have chosen to focus most of our attention in this paper on 1+1-wave maps into a certain Lorentzian 2-metric, $(L,\bsy{\lambda})$. This metric is distinguished in a number of ways.
We will show that the solution of the Cauchy problem can be expressed as the flow of a vector field which can be viewed as a curve in a certain solvable Lie algebra $\mathfrak{r}$. The Lie algebra $\mathfrak{r}$ can in turn be identified as a realisation of the Vessiot algebra associated with the nonlinear sigma model that governs wave maps into $(L,\bsy{\lambda})$. It follows that, because $\mathfrak{r}$ is solvable, the theory of systems of Lie type permits a reduction of the Cauchy problem to quadrature. This makes it possible to undertake very explicit analysis of wave maps into this Lorentzian metric. 

The plan of our paper is as follows. In section 2 we will review the main facts about Darboux integrable exterior differential systems, and in section 3 we will briefly review the known classification results on Darboux integrable surface metrics---that is, those Riemannian 2-metrics $(M,g)$ such that the nonlinear sigma models for harmonic maps $\mathbb{R}^{1,1}\to (M,g)$ are Darboux integrable at order 2. In section 4 we will study the Cauchy problem for such harmonic maps into the metric $(L,\bsy{\lambda})$; because we wish to convey a clear geometric picture about the nature of the harmonic maps we construct, we also include a subsection in which numerical solutions of an initial-boundary value problem are depicted as elastic strings vibrating on the Lorentzian surface with metric $(L,\bsy{\lambda})$. This is contrasted with the string in the Euclidean plane having precisely the same initial/boundary conditions. In section 5 we will use the theory of Darboux integrable exterior differential systems  \cite{AFV} and the generalised Goursat normal form \cite{Vassiliou2006a}, \cite{Vassiliou2006b} to construct a hyperbolic Weierstrass representation for wave maps into $(L,\bsy{\lambda})$.  Finally, section 6 is an appendix containing some technical details related to the construction in section 5.

\section{Preliminaries}
\subsection{Harmonic maps}
Let $(M,g)$ and $(N,h)$ be Riemannian or pseudo-Riemannian manifolds and $\varphi:M\to N$ a smooth map. The {\it energy} of $\varphi$ over a compact domain $\mathcal{D}\subseteq M$ is
\[
e(\varphi)=\frac{1}{2}\int_\mathcal{D}g^{ij}(x)h_{\a\b}(\varphi)\frac{\partial\varphi^\a}{\partial x_i}\frac{\partial\varphi^\b}{\partial x_j}d{\rm vol}_M.
\] 
The critical points of $e(\varphi)$ satisfy the partial differential equation (PDE)
\[
\triangle\varphi^\g+g^{ij}\L^\g_{\a\b}\frac{\partial\varphi^\a}{\partial x_i}\frac{\partial\varphi^\b}{\partial x_j}=0,
\]
where $\triangle$ is the Laplacian on $M$ and $\L^\g_{\a\b}$ the Christoffel symbols on $N$. A map $\varphi$ is said to be {\it harmonic} if it is critical for $e(\varphi)$.

Harmonic maps generalise harmonic functions and geodesics, and have been under intense study since the pioneering work of Eells and Sampson \cite{EelsSampson}; see \cite{BairdWood} and \cite{Helein} and references therein for comprehensive introductions.  
If $\rm{domain}(\varphi)=M=\mathbb{R}$, then harmonic maps are geodesic flows. If $\rm{codomain}(\varphi)=N=\mathbb{R}$, then harmonic maps are harmonic functions. 
If $(M,g)$ is pseudo-Riemannian, then harmonic maps are called {\it wave maps} or {\it nonlinear sigma models}. In this paper our focus is on wave map systems---specifically, the case {$(M,g)=(\mathbb{R}^{1,1},dx\,dy)$. Any further reference to wave maps in this paper means that the domain space $M$ is Minkowski space $\mathbb{R}^{1,1}$ with its standard flat metric.

\subsection{Darboux integrable differential systems}
As the name implies, the notion of Darboux integrability originated in the 19th century and was most significantly developed by Goursat \cite{GoursatBook}. Classically, it was a method for constructing the ``general solution" of second order PDE in one dependent and two independent variables
\[
F(x,y,u,u_x,u_y, u_{xx}, u_{xy}, u_{yy})=0
\]
that generalised the so called ``method of Monge." It relies on the notions of {\it characteristics} and {\it first integrals}. We refer the reader to \cite{GoursatBook},  \cite{IveyLandsbergBook}, \cite{Vassiliou2000}, \cite{Vassiliou2001} for further information on classical Darboux integrability. There are also extensive studies of Darboux integrable systems relevant to the equation class under study in the works \cite{Zhiber} and
\cite{ZhiberSokolov}. 

In this paper we use a new geometric formulation of Darboux integrable exterior differential systems \cite{AFV}. At the heart of this theory is the fundamental notion of a {\it Vessiot group}, which, together with systems of Lie type, are our main tools for the study of the Cauchy problem for wave maps.   

The PDE that govern wave maps can be put into the semilinear form
\[
\bsy{u}_{xy}=\bsy{f}(x,y,\bsy{u}, \bsy{u}_x, \bsy{u}_y),\ \ \ \bsy{u}, \bsy{f}\in\mathbb{R}^n
\]
via a coordinate transformation of the form $\xi = x+y$, $\tau = x-y$.
Each solution possesses a double foliation by curves called {\it  characteristics}. 
Such PDE often model wave-like phenomena,
and projection of these curves into the independent variable space describes the space-time history of the wave propagation.
The characteristics of $\bsy{u}_{xy}=\bsy{f}$ are the integral curves of a pair of rank $n+1$ distributions
\begin{equation}\label{characheristics}
H_1=\left\{D_x+D_y\bsy{f}\cdot\P {\bsy{u}_{yy}},\ \ \P {\bsy{u}_{xx}}\right\},\ \ \ H_2=\left\{D_y+D_x\bsy{f}\cdot\P {\bsy{u}_{xx}},\ \P {\bsy{u}_{yy}}\right\},
\end{equation}
where
\[
D_x=\P x+\bsy{u}_x\P {\bsy{u}}+\bsy{u}_{xx}\P {\bsy{u}_x}+\bsy{f}\P {\bsy{u}_y},\ \ \ D_y=\P y+\bsy{u}_y\P {\bsy{u}}+\bsy{f}\P {\bsy{u}_x}+\bsy{u}_{yy}\P {\bsy{u}_y}
\]
are the {\it total differential operators} along solutions of the PDE. 
Note that if $\bsy{\theta}$ is the standard Cartan codistribution for $\bsy{u}_{xy}=\bsy{f}$, then
$H_1\oplus H_2=\rm{ann}\,\bsy{\theta}$. The distributions $H_i$ are well-defined with canonical structure.

\begin{defn}
If $\Delta$ is a distribution on manifold $M$, then a function $h:M\to\mathbb{R}$ is said to be a {\it first integral} of $\Delta$ if $Xh=0$ for all $X\in \Delta$.
\end{defn}

For later use, we will briefly review the key results we require from the theory of Darboux integrable exterior differential systems, \cite{AFV}. Let $\L_i={\rm ann}\,H_i$, and denote by $H_i^\infty$ and $\L_i^\infty$ the final terms in the derived flags of $H_i, \L_i$, respectively. 

\begin{defn}
Let $\bsy{\theta}$ be a Pfaffian system on manifold $M$, and suppose that there exist auxiliary Pfaffian systems $\L_1,\L_2$ such that
\begin{enumerate}
\item[a)] $\L_1+\L_2^{\infty}=T^*M\ \ \text{and}\ \ \L_1+\L_2^{\infty}=T^*M$;
\item[b)] $\L_1^{\infty}\cap \L_2^{\infty}=\{0\}$;
\item[c)] $\L_1\cap \L_2=\bsy{\theta}$.
\end{enumerate}
Then $(M,\bsy{\theta})$ is said to be \underline{Darboux integrable} on $M$, and $(\L_1, \L_2)$ is said to be a \underline{Darboux pair} for $\bsy{\theta}$.
\end{defn}

Condition b) implies that there are no first integrals that are common to $\L_1$ and $\L_2$, while the crucial condition a) implies that each of $\L_1$ and $\L_2$ possess {\it sufficient} first integrals.    For wave map equations, we have:
 
\begin{lem}
A semilinear system $\bsy{u}_{xy}=\bsy{f}$ with $\bsy{u, f}\in \mathbb{R}^n$ is Darboux integrable at order $k$ if each of its characteristic systems $H_i$ has at least $n+1$ independent first integrals order $k$. 
\end{lem}

A key fact, proven in \cite{AFV}, is a result from the inverse problem in the theory of quotients:
\begin{thm}\label{mainAFVthm}
Let $(M,\bsy{\theta})$ be a Pfaffian system that admits a Darboux pair. Then there are Pfaffian systems $(\wh{M}_1, \wh{\bsy{\theta}}_1)$, $(\wh{M}_2, \wh{\bsy{\theta}}_2)$ which admit a common Lie group $G$ of symmetries such that:
\begin{enumerate}
\item The manifold $M$ can be locally identified as the quotient of $\wh{M}_1\times \wh{M}_2$ by a diagonal action of $G$;
\item We have the identification 
$$
\bsy{\theta}=\left(\pi_1^*\wh{\bsy{\theta}}_1+\pi_2^*\wh{\bsy{\theta}}_2\right)/G,
$$
where $\pi_i:\wh{M}_1\times \wh{M}_2\to \wh{M}_i,\ i=1,2$ are the canonical projection maps;
\item The quotient $\bsy{\pi} :\wh{M}_1\times \wh{M}_2\to M$ to the diagonal $G$-action defines a surjective superposition formula for $\bsy{\theta}$. 
\end{enumerate}
\end{thm}
\vskip 5 pt
In this context, a {\it superposition formula} for $(M,\bsy{\theta})$ is a map $\bsy{\pi} :\wh{M}_1\times \wh{M}_2\to M$ such that if $\sigma_i:\mathcal{U}_i\to \wh{M}_i$ are integral submanifolds of $(\wh{M}_i,\wh{\bsy{\theta}}_i)$, then $\bsy{\pi}\circ(\sigma_1,\sigma_2):\mathcal{U}_1\times\mathcal{U}_2\to M$ is an integral submanifold of $\bsy{\theta}$. A superposition formula is {\it surjective} if every solution of $\bsy{\theta}$ can be expressed in this form for a fixed superposition formula $\bsy{\pi}$  as $\sigma_i$ range over the integral manifolds of $\wh{\bsy{\theta}}_i$.
\vskip 5 pt
Reference \cite{AFV} is devoted to a proof of this theorem and to the identification and explicit construction of all the entities mentioned there, including the Lie group of symmetries $G$, known as the {\it Vessiot group} of the Darboux integrable exterior differential system $\bsy{\theta}$.  It is proven that the isomorphism class of the Vessiot group is an invariant of the Darboux integrable EDS $\bsy{\theta}$.

In this paper we shall construct the Vessiot group $G$ and superposition formula for the Pfaffian system $\bsy{\theta}$ governing harmonic maps $\mathbb{R}^{1,1}\to (L,\bsy{\lambda})$. We will also construct the auxilliary Pfaffian systems $(\wh{M}_i,\wh{\bsy{\theta}}_i)$. The integration of $\bsy{\theta}$ then relies upon the integration of 
$(\wh{M}_i,\wh{\bsy{\theta}}_i)$. By the generalised Goursat normal form \cite{Vassiliou2006a}, \cite{Vassiliou2006b}, we will show that for each $i=1,2$, $(\wh{M}_i,\wh{\bsy{\theta}}_i)$ is locally equivalent to a partial prolongation of the contact system on $J^1(\mathbb{R},\mathbb{R}^2)$. In turn, this will permit us to construct a hyperbolic Weierstrass type representation for $\bsy{\theta}$.

\section{Known Darboux integrable sigma models}
The first person to treat the Cauchy problem for wave maps into Riemannian targets was Gu Chao-Hao, \cite{Gu}. He established the fundamental result that for smooth initial data, wave maps into complete Riemannian metrics have long-time existence. Gu's work  initiated many further investigations where regularity constraints on the initial data have been significantly relaxed. Furthermore, some higher dimensional problems have been treated; see \cite{Struwe}.

Our interest in the present paper is in those nonlinear sigma models which are Darboux integrable. This question was recently investigated in the thesis of R. Ream \cite{Ream} who established that at order 2 or 3, the only Riemannian metrics that give rise to Darboux integrable  nonlinear sigma models are given locally by 
\begin{equation}\label{ReamMetrics}
g_1=\frac{du^2+dv^2}{1+e^{u}},\ \ \ {\rm{and}}\ \ \ \ g_2=\frac{du^2+dv^2}{1-e^{u}}.
\end{equation}
It is obvious that these metrics can be identified by a complex translation in $u$, so Ream's is a {\it real} classification. On the other hand, it is easily shown that the metrics $g_1, g_2$ are not locally isometric and that the corresponding sigma models are not locally contact equivalent. The latter follows from the fact that the corresponding Vessiot groups are not isomorphic. However, it turns out that Ream's list is not quite complete, due to an error which is easy to repair as we show below. There is an additional metric which is Darboux integrable at order 2 that is locally inequivalent to both $g_1$ and $g_2$. Before introducing it, we shall analyse the Ream metrics a little further.

The signs in front of the exponential terms in (\ref{ReamMetrics}) are isometric invariants, while the coefficients of $u$ can be rescaled without changing the Euler-Langrange equations---that is, without changing the nonlinear sigma model. We also point out that Ream's metric $g_1$ is locally equivalent to the Ricci soliton metric \cite{ChowKnopf} otherwise known as the ``Witten black hole:"
$$
g_\S=\frac{du^2+dv^2}{1+u^2+v^2}.
$$
Indeed, the transformation 
$$
u=e^{-\frac{r}{2}}\cos\frac{\theta}{2}, \ \ \ v=e^{-\frac{r}{2}}\sin\frac{\theta}{2}
$$
identifies $g_\S$ with
$$
g_3=\frac{1}{4}\frac{dr^2+d\theta^2}{1+e^{r}},
$$
which gives rise to the same nonlinear sigma model as $g_1$.

It is useful to note that, whereas $g_1$ can be represented in conformal coordinates in several ways, it also has a representation as a surface of revolution metric obtained by taking polar coordinates for $g_\S$: $x=r\cos\theta, \ y=r\sin\theta$, and then setting $r=\sinh\rho$ to transform $g_1$ to
$$
g_4=d\rho^2+(\tanh^2\rho)\, d\theta^2.
$$
Thus, the nonlinear sigma model with Langangian density
$$
L(\bsy{u})=u_xu_y+v_xv_y\tanh^2 u
$$
gives rise in spacetime coordinates $(\tau, \xi)$ to its equations of motion 

\begin{equation}
\begin{aligned}
&u_{\tau\tau}-u_{\xi\xi}=\frac{\sinh u}{\cosh^3u}(v_\tau^2-v_\xi^2),\cr
&v_{\tau\tau}-v_{\xi\xi}=\frac{2}{\cosh u\,\sinh u}(u_\tau v_\tau-u_\xi v_\xi).
\end{aligned}
\end{equation}

These equations are reminicent of a system solved in \cite{BNC}, namely, their system (4.13). However, they are not quite the same. Nevertheless, the sigma model for
$$
du^2+(\cot^2u)~dv^2
$$ 
with Lagrangian density $L=u_xu_y+v_xv_y\,\cot^2u$ has, in spacetime coordinates, the equations of motion 

\begin{equation}
\begin{aligned}
&u_{\tau\tau}-u_{\xi\xi}=-\frac{\cos u}{\sin^3u}(v_\tau^2-v_\tau^2),\cr
&v_{\tau\tau}-v_{\xi\xi}=\frac{2}{\cosh u\,\sinh u}(u_\tau v_\tau-u_\xi v_\xi).\cr
\end{aligned}
\end{equation}
These agree precisely with equations (4.13) of \cite{BNC}, and in this case the authors have constructed an explicit solution in terms of arbitary functions.   This system, too, is Darboux integrable. It can be shown that all these Darboux integrable wave map systems have Vessiot groups which are isomorphic to either $SL(2)\times\mathbb{R}$ or to $SO(3)\times\mathbb{R}$.

\section{A nonlinear sigma model with a solvable Vessiot group}

\subsection{Completing the Ream classification}
The surface of revolution metric
\begin{equation}\label{newMetric}
du^2+\frac{1}{u^2}dv^2
\end{equation}
arose in \cite{Vassiliou2008} as an application of hyperbolic reduction. It was shown that wave maps into (\ref{newMetric}) can be constructed as integrable extensions of the scalar hyperbolic PDE
\begin{equation}\label{GoursatEq}
u_{xy}=\frac{1}{u}\sqrt{1-u_x^2\strut}\sqrt{1-u_y^2},
\end{equation}
and that the Vessiot group of the integrable extension is 4-dimensional and solvable.  Equation (\ref{GoursatEq}) occurs in Goursat's well known classification of Darboux integrable scalar second order PDE in the plane. Moreover, it is straightforward to show that the nonlinear sigma model for (\ref{newMetric}),
 
\begin{equation}\label{waveEq_newMetric}
u_{xy}+\frac{v_xv_y}{u^3}=0,\ \ \ v_{xy}-\frac{\big(u_xv_y+u_yv_x\big)}{u}=0,
\end{equation}
is itself Darboux integrable \cite{Vassiliou2008}.

\begin{lem}
The metrics
$$
g_1=\frac{dx^2+dy^2}{1+e^x}, \ \ \ g_2=\frac{dx^2+dy^2}{1-e^x},\ \ \ g_P=dx^2+\frac{1}{x^2}dy^2
$$
are pairwise inequivalent. (The metrics $g_1, g_2$ can, however, be identified by the complex transformation $x\mapsto x+\pi\sqrt{-1}$.)
\end{lem}

\noindent{\it Proof.} This is a direct application of the Cartan equivalence method for 2-metrics.
Firstly, the metric $g_P=dx^2+dy^2/x^2$ has orthonormal coframe $\o^1=dx,\ \o^2=dy/x$ and lifted coframe
$$
\Theta=\{\theta^1=\cos\lambda \,\o^1-\sin\lambda\,\o^2,\ \ \theta^2=\sin\lambda\,\o^1+\cos\lambda\,\o^2\}.
$$
The connection form on the principal $SO(2)$-bundle is
$$
\o^{12}=d\lambda+\frac{1}{x}\omega^2,
$$
giving rise to the canonical structure equations
$$
\begin{aligned}
&d\theta^1=-\o^{12}\wedge\theta^2,\cr
&d\theta^2=\o^{12}\wedge\theta^1,\cr
&d\o^{12}=-\frac{2}{x^2}\theta^1\wedge\theta^2.
\end{aligned}
$$
The coframe $\{\theta^1,\ \theta^2,\ \o^{12}\}$ has rank 2 and order 1, and its invariants are spanned by the functions $K$ and $K_1$, where $K$ denotes the Gauss curvature of the metric, and the functions $K_1, K_2$ are defined by the equation $dK = K_1\, \theta^1 + K_2\, \theta^2$.  For the metric $g_P$, we have
$$
K=-\frac{2}{x^2}, \ \ K_1=\frac{4\cos\lambda}{x^3}, \ \ K_2=\frac{4\sin\lambda}{x^3},
$$
and these functions satisfy the relation
\begin{equation}
\frac{\left(\sqrt{-{K}/{2}}\right)^3}{\sqrt{K_1^2+K_2^2}}=\frac{1}{4}.
\end{equation}
Performing the equivalent calculation for metric $g_1$ of (\ref{ReamMetrics}) gives

\begin{equation}
\frac{\left(\sqrt{-{K}/{2}}\right)^3}{\sqrt{K_1^2+K_2^2}}=\frac{\sqrt{-1}\,e^{\frac{x}{2}}}{4},
\end{equation}

and for $g_2$ of (\ref{ReamMetrics}) we have
\begin{equation}
\frac{\left(\sqrt{-{K}/{2}}\right)^3}{\sqrt{K_1^2+K_2^2}}=\frac{\,e^{\frac{x}{2}}}{4}.
\end{equation}

Since these functional relations are different, it follows from Cartan's solution of the local equivalence problem that the metrics are pairwise inequivalent. \hfill\qed

\begin{cor}
The nonlinear sigma model (\ref{waveEq_newMetric}) for wave maps into the metric
$$
g_P=du^2+\frac{1}{u^2}dv^2
$$
extends the Ream classification \cite{Ream} of metrics for which the wave map equation is Darboux integrable at order 2 (or 3).
\end{cor}

\begin{rem}
The author in \cite{Ream} missed the metric $g_P$ in his classification because of an allowable form of an integating factor with a free parameter $C$ in which he inadvertently ommitted the value $C=0$ in the prelude to equation (4.55) on p53. Setting $C=0$ prior to (4.55), following his subsequent calculations and then making an elementary change of variable leads to metric $g_P$.
\end{rem}

\subsection{The Cauchy problem}

According to Gu's theorem \cite{Gu}, the Cauchy problem for wave maps into complete Riemannian metrics has global existence for smooth initial data. The situation appears to be rather different for wave maps into split signature metrics. For instance, Terng and Uhlenbeck discovered a  bi-invariant Lorentzian metric on $SL(2,\mathbb{R})$ and showed that there were smooth, finite energy wave maps which blew up in finite time.

With the ultimate aim of studying the properties of harmonic maps between pseudo-Riemannian manifolds, we will in the remainder of this paper turn our attention to wave maps into Lorentzian metrics.  Metric (\ref{newMetric}) can be put into the conformal form
\begin{equation}\label{newMetricConformal}
\frac{d\bar{u}^2+d\bar{v}^2}{2\bar{u}}
\end{equation}
via the transformation
\begin{equation}\label{conformalTrans}
u=\sqrt{2\bar{u}},\ v=\bar{v}.
\end{equation}
In this section we will show that the general Cauchy problem for all harmonic maps from Minkowski space $\mathbb{R}^{1,1}$ into the
semi-Riemannian metric 
\begin{equation}\label{semiRiemMetricConformal}
\frac{du_1^2-du_2^2}{2u_1}.
\end{equation}
is reducible to quadrature and that, given any initial data, its solution can be constructed by solving the initial value problem for an ordinary differential equation of {\it Lie type} \cite{Bryant}, \cite{Carinena}, \cite{DoubrovKomrakov}, \cite{Vassiliou2012} associated to the free and transitive action of a solvable Lie group on $\mathbb{R}^4$.
\vskip 10 pt

The nonlinear sigma model in this case is also Darboux integrable. In lightcone coordinates on both source $(x,y)$ and target $(u,v)$ spaces
$$
(u_1,u_2)\mapsto\left(\frac{u_1+u_2}{2},\ \frac{u_1-u_2}{2}\right)=(u,v)
$$
it has the form
\begin{equation}\label{mainEquations}
u_{xy}=\frac{u_xu_y}{u+v},\ \ \ \ v_{xy}=\frac{v_xv_y}{u+v}.
\end{equation}

Each of the characteristic systems $H_1$ and $H_2$ in  (\ref{characheristics}) have 4 independent first integrals, which we shall label
$$
y,\ \beta_1,\ \beta_2,\ \beta_3\ \ {\rm for} \ \ H_1\ \ \ {\rm and}\ \ \ x,\ \alpha_1,\ \alpha_2,\ \alpha_3\ \ {\rm for} \ \ H_2.
$$
Specifically, we have
$$
\beta_1=\frac{u_yv_y}{u+v},\ \ \ \ \beta_2=\frac{d\beta_1}{dy},\ \ \ \ \beta_3=\frac{v_{yy}}{v_y}-\frac{u_y}{u+v},
$$
and
$$
\alpha_1=\frac{u_xv_x}{u+v},\ \ \ \ \alpha_2=\frac{d\alpha_1}{dx},\ \ \ \ \alpha_3=\frac{u_{xx}}{u_x}-\frac{v_x}{u+v}.
$$
Adapting these first integrals as new coordinates on the differential equation submanifold of $J^2(\mathbb{R}^2,\mathbb{R}^2)$ considerably simplifies the Pfaffian system 
for (\ref{mainEquations}). We obtain:

\begin{eqnarray}
\bsy{\theta}=\Bigg\{db_1-b_2dy,\ \ \ \  da_1-a_2dx,\ \  \ \ \ dz_1-z_3dx-\frac{b_1}{z_4}(z_1+z_2)dy,\cr
 dz_2-\frac{a_1}{z_3}(z_1+z_2)dx-z_4dy, \ \ \  dz_3-(a_1+z_3a_3)dx-\frac{b_1z_3}{z_4}dy,\cr
 \ \ \ \ \  dz_4-\frac{a_1z_4}{z_3}dx-(b_3z_4+b_1)dy\Bigg\},\cr
\end{eqnarray}
where
$$
z_1=u,\ z_2=v,\ z_3=u_x,\ z_4=v_y,\ a_i=\alpha_i,\ b_j=\beta_j. 
$$

Our main results are contained in the following theorem:

\begin{thm}\label{IVPmain}
Consider the initial value problem for the harmonic map system
$$
u_{xy}=\frac{u_xu_y}{u+v},\ \ \ \ v_{xy}=\frac{v_xv_y}{u+v},
$$
with initial data
\begin{equation}\label{IC}
u_{|_\g}=\phi_1,\ \  v_{|_\g}=\phi_2,\ \ \frac{\partial u}{\partial\boldsymbol{n}}_{|_\g}=\psi_1,\ \ \frac{\partial v}{\partial\boldsymbol{n}}_{|_\g}=\psi_2,
\end{equation}
where $\g$ is a curve with tangents nowhere parallel to the $x$- or $y$-axes, $\bsy{n}$ is a unit normal vector field along $\g$ and $\phi_i, \psi_i$ are smooth functions along $\g$. 

\begin{enumerate}
\item The Cauchy problem  has a unique smooth local solution. Moreover, the unique local solution is expressible as the solution of an ordinary differential equation of Lie type associated to a local action of a solvable Lie group on $\mathbb{R}^4$. 

\item Given the unique local solution $(u,v)$ from part 1., the Cauchy problem for harmonic maps
\[
(\mathbb{R}^{1,1},\ dxdy)\to \left(M, \frac{du_1^2-du_2^2}{2u_1}\right)
\]
is given by
\[
u_1=u+v,\ \ u_2=u-v,
\]
where $u_1,u_2$ satisfy initial conditions
\[
{u_1}_{|_\g}=\phi_1+\phi_2,\ \  {u_2}_{|_\g}=\phi_1-\phi_2,\ \ \frac{\partial u_1}{\partial\boldsymbol{n}}_{|_\g}=\psi_1+\psi_2,\ \ \frac{\partial u_2}{\partial\boldsymbol{n}}_{|_\g}=\psi_1-\psi_2.
\]
\end{enumerate}

\end{thm}

\noindent{\it Proof.} By hyperbolicity, the problem is locally well posed. Let $k_1(y), k_2(y)$ be arbitrary, smooth real-valued functions, and consider the overdetermined PDE system defined by (\ref{mainEquations}) together with the additional PDE
\begin{equation}\label{Cau2}
\b_1=k_1(y),\ \b_2=\dot{k}_1(y),\ \b_3=k_2(y),
\end{equation}
where the dot denotes $y$-differentiation. It can be shown that {\it this} overdetermined system $\mathcal{E}'$ is involutive and, moreover, admits a 1-dimensional Cauchy characteristic distribution. Now suppose we fix a smooth curve $\g$ embedded in a portion of the $xy$-plane $\mathcal{N}$, and suppose that Cauchy data is prescribed along $\g$ as described in (\ref{IC}). Then by an argument similar to (\cite{Vassiliou2000}, Proposition 3.3), $\g$ can be lifted to a unique curve $\wh{\g} : I\subseteq\mathbb{R}\to J^1(\mathcal{N},\mathbb{R}^2)$ which agrees with the Cauchy data. Let $\iota : \mathcal{H}_{k_1,k_2}\to J^2(\mathcal{N},\mathbb{R}^2)$ denote the submanifold in $J^2(\mathcal{N},\mathbb{R}^2)$ defined by PDE system $\mathcal{E}'$, and let $\Theta$ denote the contact system on $J^2$. Let
$\bsy{\theta}_{k_1,k_2}=\iota^*\Theta$ be the Pfaffian system whose integral submanifolds are the solutions of $\mathcal{E}'$. We now claim that the functions $k_1, k_2$ can be uniquely chosen in such a way that  $\wh{\g}$ can be extended to a 1-dimensional integral submanifold $\widetilde{\g}$ of $\bsy{\theta}_{k_1,k_2}$. To see this, consider the Pfaffian system
\begin{eqnarray}
\bsy{\theta}_{k_1,k_2}=\Big\{ \o^1=du-u_xdx-u_ydy,\ \o^2=dv-v_xdx-v_ydy,\cr
\ \ \o^3=du_x-u_{xx}dx-\frac{u_x u_y}{u+v}  dy,\ \ \ \o^4=du_y- \frac{u_x u_y}{u+v}dx-u_{yy}dy,\\[0.1in]
\ \ \o^5=dv_x-v_{xx}dx- \frac{v_x v_y}{u+v}dy,\ \  \o^6=dv_y- \frac{v_x v_y}{u+v}dx-v_{yy}dy\Big\}.\cr
\end{eqnarray}
Pulling this system back by $\wh{\g}$, we observe that $\o^1, \o^2$ pull back to zero by construction. The forms $\o^4$ and $\o^6$ define the functions $u_{yy}$ and $v_{yy}$ along $\g$,
while $\o^3, \o^5$ define the functions $u_{xx}$ and $v_{xx}$ along $\g$. All these functions are expressed in terms of the Cauchy data $\phi_i, \psi_i$. Substituting these expression back into (\ref{Cau2}) uniquely determines the functions $k_1, k_2$ in terms of $\phi_i, \psi_i$.

We will now use the first integrals $\beta_j$ recorded above to demonstrate that $\bsy{\theta}_{k_1,k_2}$ has a one-dimensional Cauchy characteristic distribution and that, in particular, the Cauchy characteristic vector field can be chosen to be a curve in a certain Lie algebra---the Vessiot algebra \cite{AFV} of system (\ref{mainEquations}). It will be seen that the Cauchy characteristic vector field is generically transverse to the Cauchy data and extends the one-dimensional integral submanifold of $\bsy{\theta}_{k_1,k_2}$ to the solution of the Cauchy problem. Because the Cauchy characteristic vector field is a curve in a Lie algebra, this extension from a one-dimensional to a two-dimensional integral of $\bsy{\theta}_{k_1,k_2}$ is an ordinary differential equation $\mathfrak{L}$ of Lie type. Its coefficients and initial conditions are fixed by all the data present in the problem, including the Cauchy data. {\it Any} solution of $\mathfrak{L}$ (independently of its initial conditions) permits a Lie reduction of $\mathfrak{L}$ and will permit us to solve the IVP for $\mathfrak{L}$.

We find by explicit calculation that the vector field

\begin{equation}\label{CauchyVector}
\xi_{k_1,k_2}=\P y-R_1-k_2(y)R_2+k_1(y)R_3
\end{equation}
spans the Cauchy characteristic vectors of $\bsy{\theta}_{k_1,k_2}$, where $\mathfrak{r}=\{R_1, R_2, R_3, R_4\}$ is a basis for a solvable Lie algebra with nonzero structure
$$
[R_1,R_2]=R_1,\ \ [R_1,R_3]=-R_4,\ \ [R_2,R_3]=R_3
$$
with local expression
\begin{equation}
R_1=-z_4\P {z_2}, R_2=-z_4\P {z_4},\ R_3=\frac{z_1+z_2}{z_4}\P {z_1}+\frac{z_3}{z_4}\P {z_3}+\P {z_4},R_4=\P {z_1}-\P {z_2}.
\end{equation}
Note that $\mathfrak{r}$ generates a free transitive local Lie group action on $\mathbb{R}^4$. 

We have shown that, given Cauchy data $\phi_i, \psi_i$ along a generic smooth curve $\g$ in the $xy$-plane, there is a unique choice of $k_1$ and $k_2$ such that the Cauchy data can be lifted to a 1-dimensional integral submanifold of $\bsy{\theta}_{k_1,k_2}$. Indeed, following the outlined procedure in the case of the important curve $\g(x)=(x,x)$ gives
$$
k_1(x)=\frac{2\psi_1\psi_2+{\phi_1}_x{\phi_2}_x-\sqrt{2}({\phi_1}_x\psi_2+\psi_1{\phi_2}_x)}{4(\phi_1+\phi_2)}
$$
and
\[ 
k_2(x) = \frac{(\phi_1 + \phi_2)(2\phi_{2xx} - 2 \sqrt{2} \psi_{2x}) + (\sqrt{2}(\psi_1 - \psi_2) - (\phi_{1x} + \phi_{2x}))(\phi_{2x} - \sqrt{2} \psi_2)  }{\sqrt{2}(\phi_1 + \phi_2)(\sqrt{2} \phi_{2x} - 2\psi_2)} .
\]

Thus for any smooth, real-valued functions $k_1(y), k_2(y)$, $\bsy{\theta}_{k_1,k_2}$ admits a Cauchy characteristic vector field $\xi_{k_1,k_2}$. Furthermore, we have shown that the Cauchy data along $\gamma$ can always be lifted to a unique 1-dimensional integral manifold of $\bsy{\theta}_{k_1,k_2}$ for a unique choice of $k_1,k_2$, depending upon the Cauchy data. The Cauchy characteristic vector field can then be used to extend this $1$-dimensional integral manifold to a $2$-dimensional integral manifold, provided that the Cauchy characteristic vector field is everywhere transverse to the $1$-dimensional manifold.
But because we have chosen initial data which has tangents parallel to neither coordinate axis in the Minkowski plane, transversality is guaranteed; hence, the Cauchy characteristic vector field $\xi_{k_1,k_2}$ provides an extension of the initial data integral manifold $\widetilde{\g}$ to the unique smooth local solution of the Cauchy problem.

To see that the extension by $\xi_{k_1,k_2}$ from a 1- to a 2-dimensional integral manifold is an ODE of Lie type, it is sufficient to note that $\xi_{k_1,k_2}$ is a curve in the Lie algebra $\mathfrak{r}$. \hfill\qed

\begin{rem}
At this point, the occurrence of Lie algebra $\mathfrak{r}$ may appear somewhat mysterious. Its geometric origin and construction is discussed in the Appendix in section 6 and used again in section 5.
\end{rem}

\begin{lem}
The metric $$\bsy{\lambda}=\frac{du_1^2-du_2^2}{2u_1}$$ is complete.
\end{lem}

\noindent{\it Proof.} In lightcone coordinates $(u,v)$ on $(L, \bsy{\lambda})$, the geodesic equations are
$$
\ddot{u}=\frac{\dot{u}^2}{u+v},\ \  \ddot{v}=\frac{\dot{v}^2}{u+v}.
$$
The complete set of geodesics are
$$
u=\frac{1}{(a-b)^2}\bigg(a^2(a_1+b_1)e^{(a-b)t}+ab(b-a)(a_1+b_1)t-a^2b_1-2aba_1+a_1b^2\bigg),
$$
$$
v=\frac{1}{(a-b)^2}\bigg(b^2(a_1+b_1)e^{-(a-b)t}-ab(b-a)(a_1+b_1)t+a^2b_1-2abb_1-a_1b^2\bigg),
$$
where $a, b, a_1, b_1$ are constants depending upon initial conditions. Note that $u={\rm constant}$ and $v={\rm constant}$ are geodesics of $\bsy{\lambda}$.
\hfill\qed
\vskip 10 pt
\noindent{\bf Example 4.1.}
We will use the above solution of the Cauchy problem to show that, even with real analytic initial data, harmonic maps from Minkowski space into the Lorentzian metric $\bsy{\lambda}$ may develop  shocks  in finite time, despite the fact that the metric is complete, in contrast to Gu's theorem \cite{Gu}.

\vskip 10 pt

Consider the initial value problem  for (\ref{mainEquations}) subject to the initial conditions
\begin{equation}\label{IVP1}
u_{|_\g}=1-x,\ \  v_{|_\g}=2x,\ \ \frac{\partial u}{\partial\boldsymbol{n}}_{|_\g}=0,\ \ \frac{\partial v}{\partial\boldsymbol{n}}_{|_\g}=0,
\end{equation}
where $\g$ is the curve $y=x$ in the $xy$-plane and $\bsy{n}=2^{-\frac{1}{2}}(\P x-\P y)$.  We lift $\g$ to a 1-dimensional integral of $\bsy{\theta}_{k_1,k_2}$ for suitable unique functions $k_1, k_2$, as follows: observe that
$$
{z_3}_{|_\g}={u_x}_{|_\g}=\frac{1}{2}{\phi_1}_x+\frac{1}{\sqrt{2}}\psi_1,\ {z_4}_{|_\g}={v_y}_{|_\g}=\frac{1}{2}{\phi_2}_x+\frac{1}{\sqrt{2}}\psi_2 .
$$
Hence we have the partial lift
$$
(z_1, z_2, z_3, z_4)_{|_\g}=\left(1-x,\ 2x,\ -\frac{1}{2},\ 1\right). 
$$
We can extend this partial 1-dimensional integral by direct substitution, and we find that the curve $\sigma:\mathbb{R}\to \mathcal{H}_{k_1,k_2}$ defined by
$$
\begin{aligned}
(z_1, z_2, z_3, z_4, a_1, a_2, a_3, b_1, b_2, &b_3)_{|_\g}=\Bigg(1-x, 2x, -\frac{1}{2}, 1, -\frac{1}{2+2x}, \frac{2}{(2+2x)^2},\cr 
 &-\frac{1}{2+2x}, -\frac{1}{2+2x}, \frac{2}{(2+2x)^2}, -\frac{1}{2+2x}\Bigg)\cr
\end{aligned}
$$
along $y=x$ is a 1-dimensional integral manifold of $\bsy{\theta}_{k_1,k_2}$ where
$$
k_1(y)=k_2(y)=-\frac{1}{2+2y}.
$$
Hence the Cauchy characteristic vector field of  $\bsy{\theta}_{k_1,k_2}$ is 
\begin{equation}\label{CauchyVecEx1}
\P y-R_1+\frac{1}{2+2y}R_2-\frac{1}{2+2y}R_3.
\end{equation}
Computing the flow for (\ref{CauchyVecEx1}) subject to the initial conditions above along $\g$ gives the explicit solution of Cauchy problem (\ref{IVP1}):
\begin{eqnarray}
u &=& -\frac{(-y\sqrt{1+x}-\sqrt{1+y}+x\sqrt{1+y}-\sqrt{1+x}+(1+y)^{3/2}}{\sqrt{1+y}},\nonumber\cr \cr \cr
v &=& -\frac{4\sqrt{1+y}-4\sqrt{1+x}+y\sqrt{1+y}-4y\sqrt{1+x}+x\sqrt{1+y}}{\sqrt{1+y}}.\cr
\end{eqnarray}
The corresponding solution $u_1=u+v, u_2=u-v$ for the metric (\ref{semiRiemMetricConformal}) in space-time coordinates
$$
x=\frac{\xi+\tau}{2},\ \ y=\frac{\xi-\tau}{2}
$$
is
\begin{eqnarray}
u_1&=&\frac{(5\xi-5\tau+10)\sqrt{4+2\xi+2\tau}-(8+4\xi)\sqrt{4+2\xi-2\tau}}{2\sqrt{4+2\xi-2\tau}}, \cr
\cr \cr 
u_2&=&\frac{(3\tau-3\xi-6)\sqrt{4+2\xi+2\tau}+8\sqrt{4+2\xi-2\tau}}{2\sqrt{4+2\xi-2\tau}} .\nonumber\cr
\end{eqnarray}

That is, the functions $u_i(\xi,\tau)$ satisfy
\begin{equation}\label{full_spacetime}
\begin{aligned}
&\frac{\partial^2u_1}{\partial\tau^2}-\frac{\partial^2u_1}{\partial\xi^2}=\frac{1}{2u_1}\left(\left(\frac{\partial u_1}{\partial\tau}\right)^2+
\left(\frac{\partial u_2}{\partial\tau}\right)^2-\left(\frac{\partial u_1}{\partial\xi}\right)^2-\left(\frac{\partial u_2}{\partial\xi}\right)^2\right),\cr
&\frac{\partial^2u_2}{\partial\tau^2}-\frac{\partial^2u_2}{\partial\xi^2}=
\frac{1}{u_1}\left(\frac{\partial u_1}{\partial\tau}\frac{\partial u_2}{\partial\tau}-
\frac{\partial u_1}{\partial\xi}\frac{\partial u_2}{\partial\xi}\right).
\end{aligned}
\end{equation} 

Note that in space-time coordinates, the initial conditions have the form
$$
u_1(\xi,0)=1+\frac{1}{2}\xi,\ \ u_2(\xi,0)=1-\frac{3}{2}\xi,\ \ \frac{\partial u_1}{\partial\tau}(\xi,0)=\frac{\partial u_2}{\partial\tau}(\xi,0)=0.
$$

Despite the analyticity of this initial data, the solution exhibits a shock wave whose image is moving in the $\xi\cdot\tau$-plane along the curve 
$$
\tau=\frac{1}{2}\sqrt{8\xi+3\xi^2}.
$$  
At a given value of $\xi>0$, both solutions $u_1, u_2$ experience infinite discontinuities at a time $\tau$ soon after $\tau=\frac{1}{2}\sqrt{8\xi+3\xi^2}$. This behaviour is excluded in the case of harmonic maps from Minkowski space to complete Riemannian metrics according to the result of Gu \cite{Gu}.

\vskip 10 pt
\noindent{\bf Example 4.2.}
On the other hand, if we take a slight variation of these initial conditions and choose
\begin{equation}\label{IVP2}
u_{|_\g}=1-x^2,\ \  v_{|_\g}=x^2,\ \ \frac{\partial u}{\partial\boldsymbol{n}}_{|_\g}=0,\ \ \frac{\partial v}{\partial\boldsymbol{n}}_{|_\g}=0,
\end{equation}
we find that
\[ k_1(y) = -y^2, \ \ \ k_2(y) = \frac{1}{y}. \]
Then computations similar to those of the previous example show that the Cauchy characteristic vector field \begin{equation}\label{CauchyVecEx2}
\P y-R_1-\frac{1}{y}R_2-y^2R_3
\end{equation}
of $\bsy{\theta}_{-y^2,1/y}$ leads to the simple globally defined solutions:
\begin{equation}
u_1=1-\frac{1}{4}\xi^2\tau^2,\ \ u_2=1-\frac{1}{2}(\xi^2+\tau^2)
\end{equation}
satisfying the following initial conditions in space-time coordinates:
$$
u_1(\xi,0)=1,\ \ u_2(\xi,0)=1-\frac{1}{2}\xi^2,\ \ \frac{\partial u_1}{\partial\tau}(\xi,0)=\frac{\partial u_2}{\partial\tau}(\xi,0)=0.
$$

Note that in both Examples 4.1 and 4.2, the energy of the initial data is infinite. However, in the latter case the energy density is an analytic function of the parameter $\xi$, whereas in the former it has a simple pole in $\xi$.

\subsection{The Cauchy problem and systems of Lie type}

One of the reasons we are able to get compact, simple formula for solutions of the Cauchy problem is that the Vessiot group $G$ (whose Vessiot algebra is $\mathfrak{r}$) is solvable. As we have shown, the solution of the Cauchy problem is equivalent to the flow of a specially constructed vector field $\xi_{k_1,k_2}$ that depends on the choice of initial conditions 
$\phi_i, \psi_i$ and their 2-jets. We have observed that $\xi_{k_1,k_2}$ is a curve in the Vessiot algebra $\mathfrak{r}$ of the nonlinear sigma model (\ref{mainEquations}).  That is, the differential system that explicitly determines the solution of the general Cauchy problem is an ODE of Lie type. It is well known 
\cite{DoubrovKomrakov}, \cite{Bryant}, \cite{Vassiliou2012} that if, as in this case, the group action associated to the system of Lie type is solvable, then the solution can be reduced  to quadrature. Thus, we have

\begin{prop}
The solution of the Cauchy problem for wave maps of Theorem \ref{IVPmain} 
\[
(\mathbb{R}^{1,1},\ dx\,dy)\to \left(M, \frac{du_1^2-du_2^2}{2u_1}\right)
\]
is reducible to quadrature. The solution can constructed by extending the Cauchy data to a 1-dimensional integral manifold of $\bsy{\theta}_{k_1,k_2}$ and then flowing this out to the unique solution via the flow of the Cauchy vector $\xi_{k_1,k_2}$ (of $\bsy{\theta}_{k_1,k_2}$), where $k_1(y),k_2(y)$ are real-valued functions that are uniquely determined by the initial data. 
\end{prop}

\subsection{Numerical example}

Recall from the Introduction that an elastic string vibrating without friction and constrained on a surface equipped with a Riemannian or pseudo-Riemannian metric $g$ is mathematically modelled by wave maps into $g$. To gain additional insight into the nonlinear dynamical systems under study, we include a brief section on the interpretation of nonlinear sigma models as partial differential equations that {\it precisely} describe this motion. We are interested in particular in the nonlinear sigma model for harmonic maps from Minkowski space $\mathbb{R}^{1,1}$ into the metric $\bsy{\lambda}$. In order that boundary conditions also be included, we have computed and graphed numerical solutions of two wave map systems. Figures 1 and 2 each display graphs of wave maps subject to the boundary and initial conditions
$$
\begin{aligned}
&u_1(0,\tau)=0,\ u_1(2,\tau)=2,\ u_2(0,\tau)=1,\ u_2(2,\tau)=1,\cr 
& u_1(\xi,0)=\xi,\ u_2(\xi,0)=1,\ \frac{\partial u_1}{\partial\tau}(\xi,0)=0,\ \ \frac{\partial u_2}{\partial\tau}(\xi,0)=\xi(\xi-2)
\end{aligned}
$$ 
 
\begin{tabular}{ccc}
\includegraphics[width=1.8in]{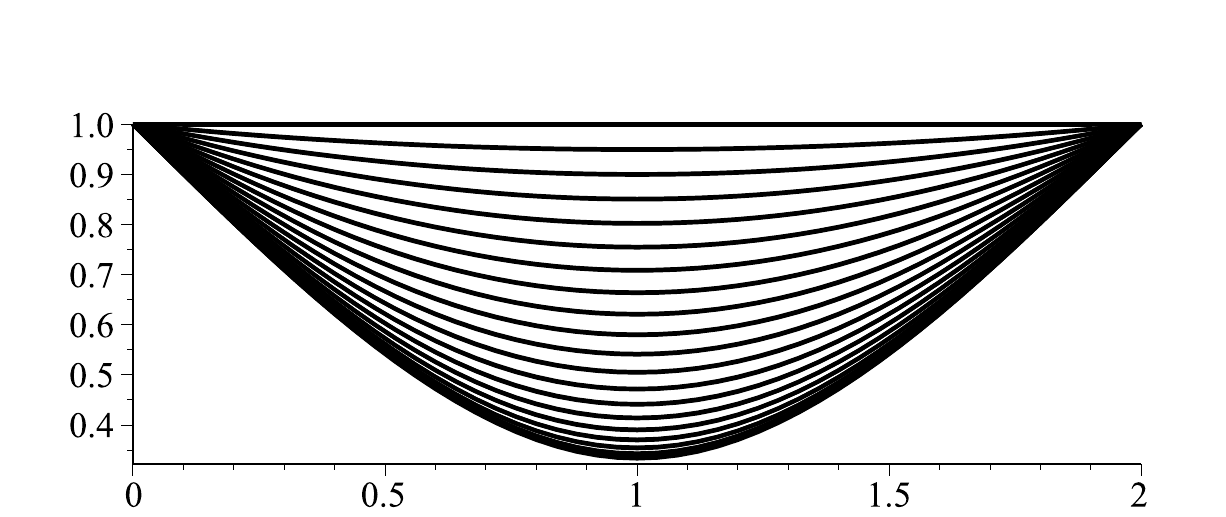}&\includegraphics[width=1.8in]{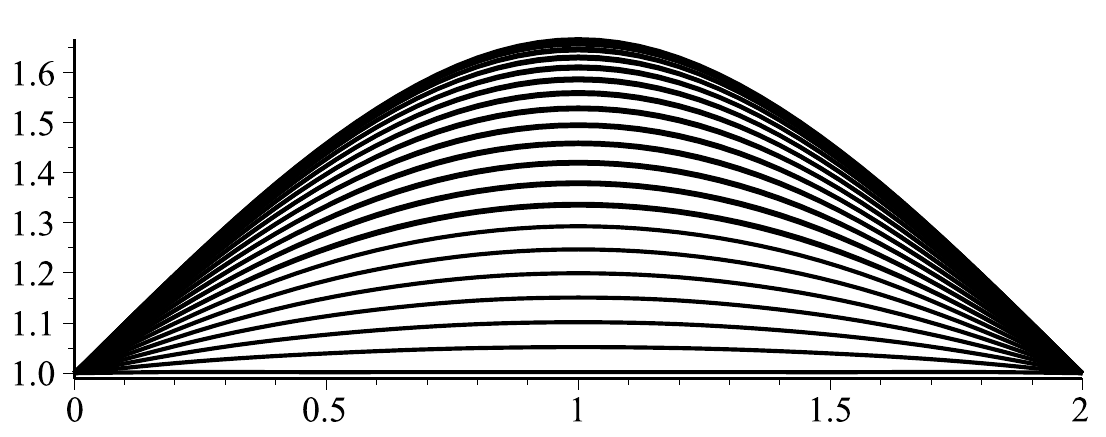}&\includegraphics[width=1.8in]{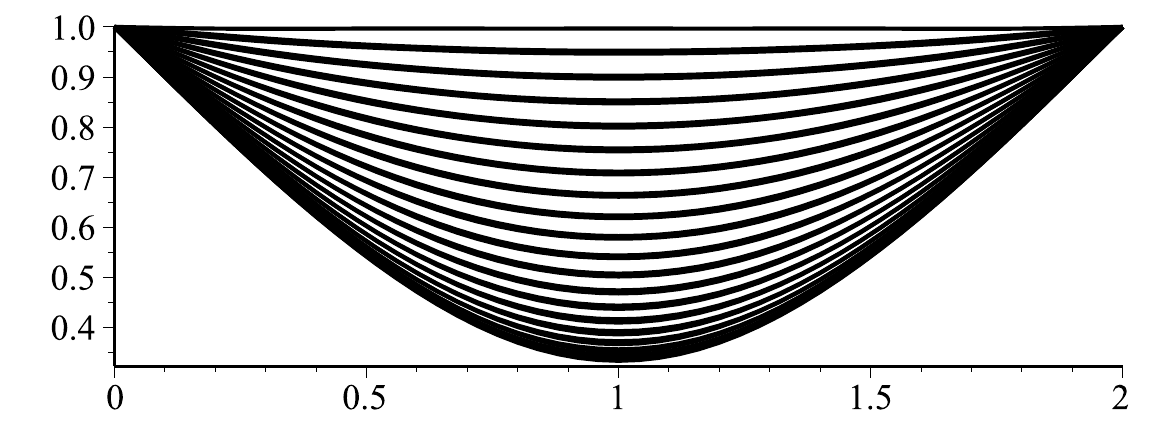}\cr
\includegraphics[width=1.8in]{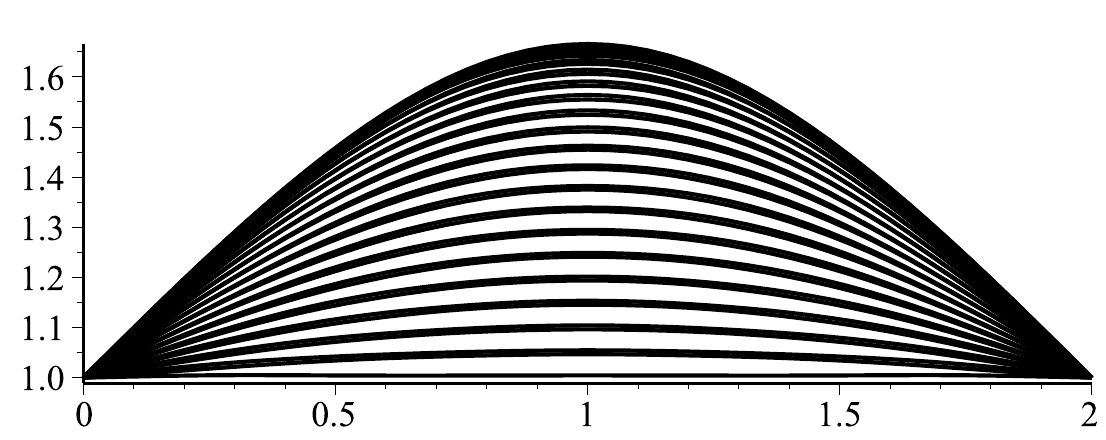}&\includegraphics[width=1.8in]{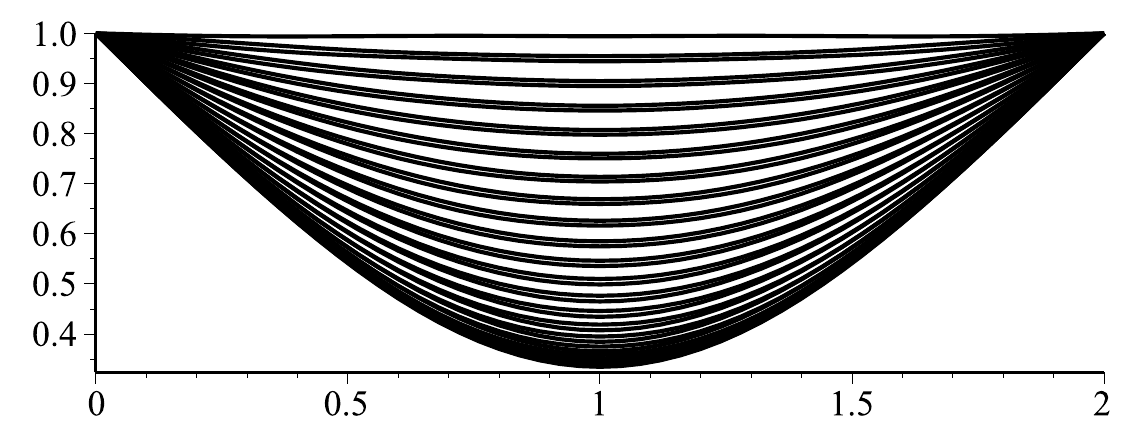}&\phantom{jgjhgjhgjhgjh}\cr
\end{tabular}
\begin{center}
Figure 1: Evolution of the string in Euclidean plane.
\end{center}

\begin{tabular}{ccc}
\includegraphics[width=1.8in]{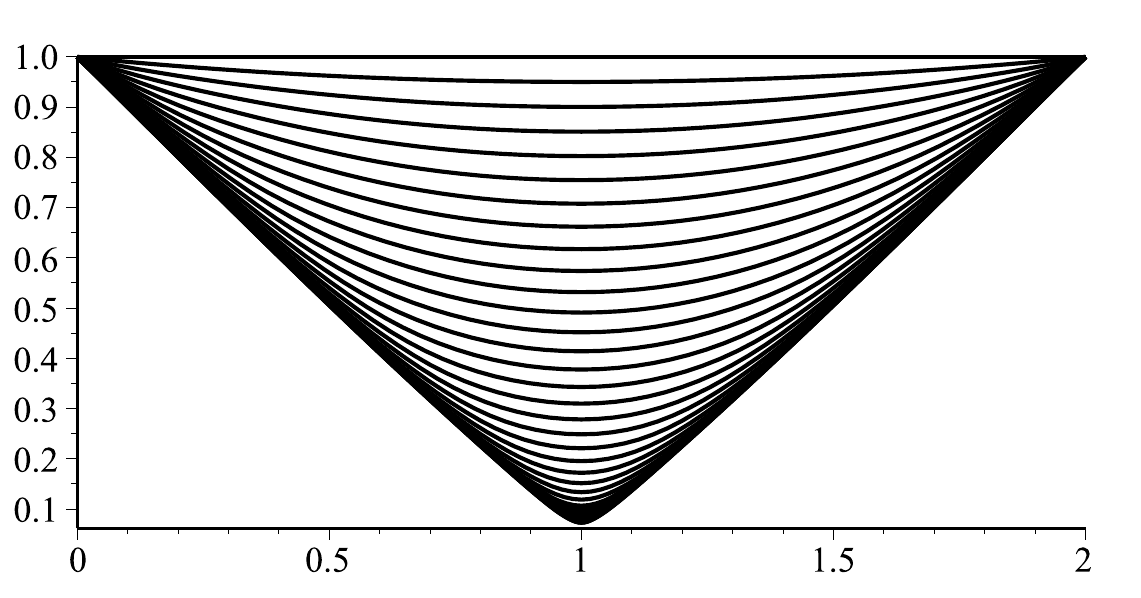}&\includegraphics[width=1.8in]{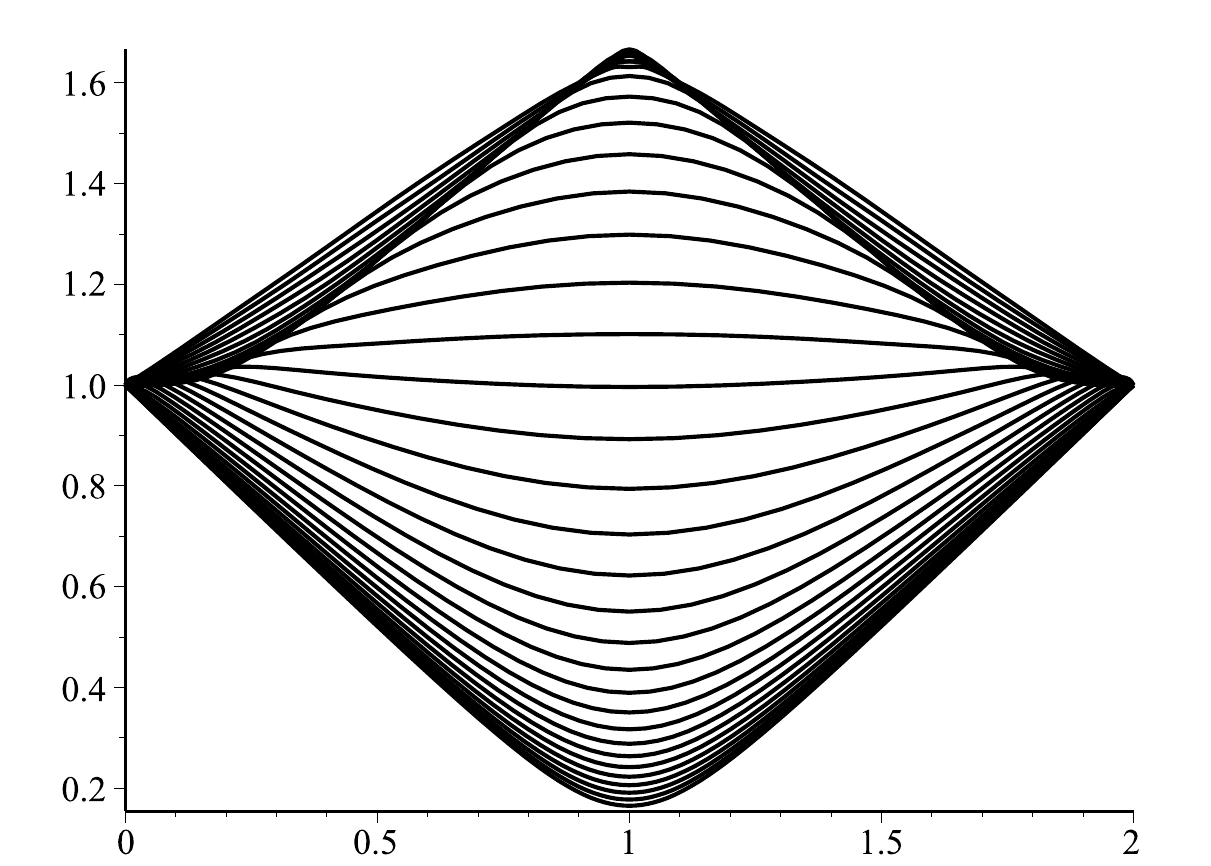}&
\includegraphics[width=1.8in]{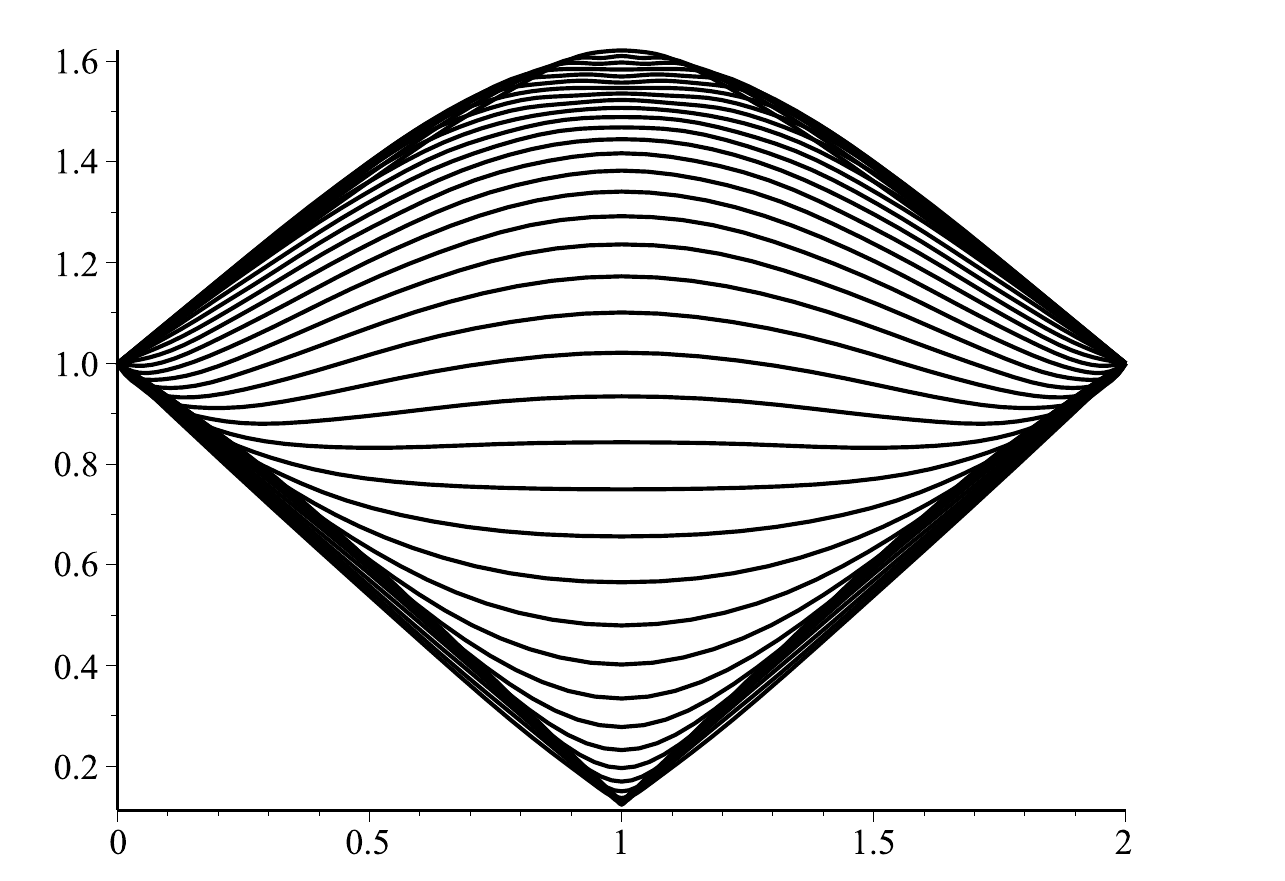}\cr
\includegraphics[width=1.8in]{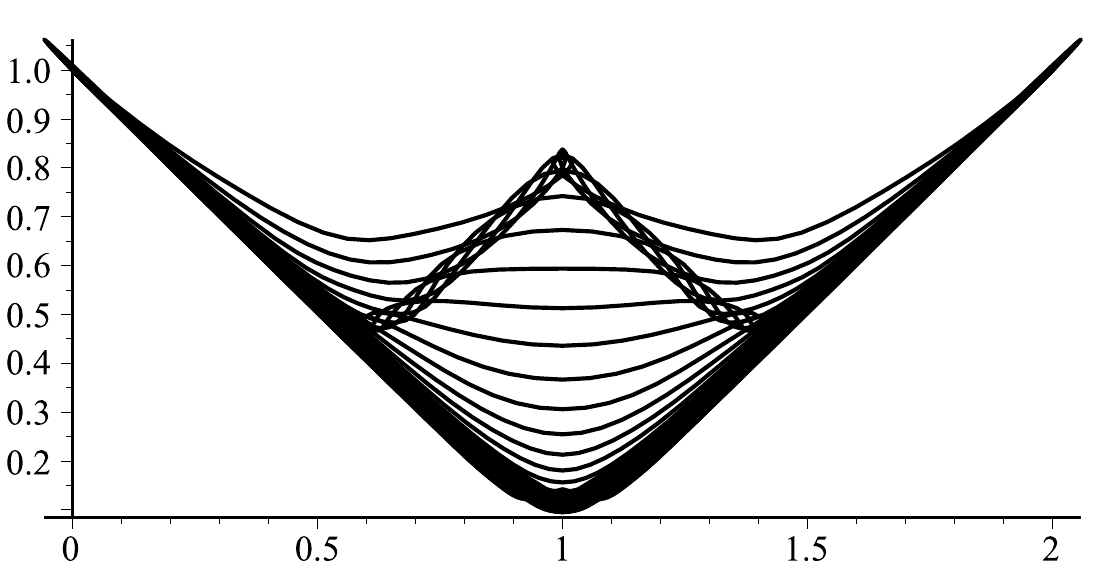}&\includegraphics[width=1.8in]{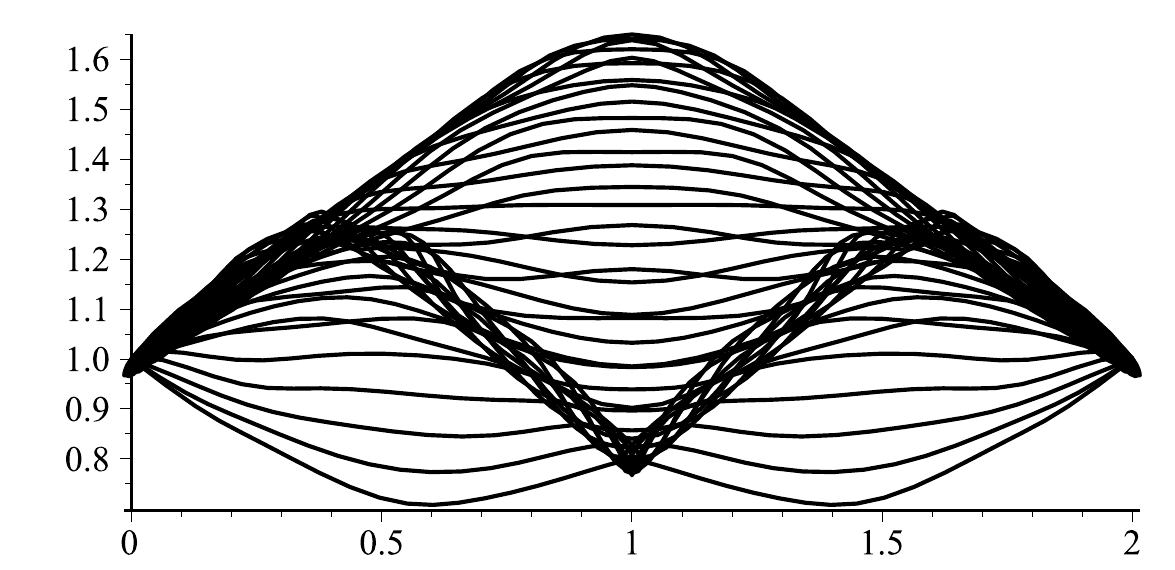}&\phantom{jhkjhggjkhgkjhg}\cr
\end{tabular}
\begin{center}
Figure 2: Evolution of the string in Lorentzian metric $\bsy{\lambda}$.
\end{center}


The values of $u_2$ are plotted along the vertical axis and those of $u_1$ along the horizontal axis. Initially, the string occupies 
$0\leq u_1\leq 2, u_2=1$; is clamped at $(u_1,u_2)=(0,1)$ and $(u_1,u_2)=(2,1)$, and the $u_2$-component has initial velocity $u_1(u_1-2)$ at $u_1$.
Figure 2 plots the time $(\tau)$ evolution of curves $\bsy{u}: \mathbb{R}^{1,1}_{\xi,\tau}\to (L, \bsy{\lambda})$ for 
$\tau\in [0,T]$, for some $T>0$. For comparison, exactly the same boundary and initial conditions are imposed on wave maps into the {\it Euclidean plane} and displayed in Figure 1. Note that $u_2=1$ is not a geodesic of $\bsy{\lambda}$, so even without the initial downward initial velocity the string would vibrate from its initial position. However, the nonzero initial velocity increases the amplitude of the vibration and leads to more interesting behaviour. By the fourth frame in Figure 2, the curve/string ceases to be immersed at $u_1=1$ and at some epochs during its motion. It also loops around to intersect itself briefly. At this time, symmetric ``waves" form on each side which persist to the end of the computation. The lifetime of this solution is not known to the authors.\footnote{See \cite{ClellandVassiliou} for further discussion and examples of string motion on Riemannian surfaces.}

\section{Hyperbolic Weierstrass representation}

In this section we use the Darboux integrability of the wave map equation (\ref{mainEquations}) to compute its general solution and hence construct a hyperbolic Weierstrass-type representation
for wave maps into the corresponding metric. According to \cite{AFV}, we pull back $\bsy{\theta}$ to suitable integral manifolds $M_1, M_2$ of $H_1^{(\infty)}$ and $H_2^{(\infty)}$ 
respectively. It is convenient to define $M_1$ by $y=b_1=b_2=b_3=0$ and $M_2$ by $x=a_1=a_2=a_3=0$. This gives Pfaffian systems
$$
\begin{aligned}
\bsy{\theta}_1=
\Bigg\{da_1-a_2dx, dz_1-z_3dx, dz_2-\frac{a_1}{z_3}(z_1+z_2)dx, dz_3-(a_1&+a_3z_3)dx,\cr 
& dz_4-\frac{a_1z_4}{z_3}dx\Bigg\}\cr
\end{aligned}
$$
and
$$
\begin{aligned}
\bsy{\theta}_2=
\Bigg\{db_1-b_2dy, dz_1-\frac{b_1}{z_4}(z_1+z_2)dy, dz_2-z_4dy, dz_3-&\frac{b_1z_3}{z_4}dy,\cr 
&dz_4-(b_1+b_3z_4)dy\Bigg\},\cr
\end{aligned}
$$
each of rank 5, on 8-manifolds $(M_i,\ \bsy{\theta}_i)$. Locally $M_1$ has coordinates 
$$
x, z_1, z_1, z_2, z_3, z_4, a_1, a_2, a_3,
$$
while $M_2$ has local coordinates
$$
y, z_1, z_1, z_2, z_3, z_4, b_1, b_2, b_3.
$$
Using these local formulas, we define a local product structure 
\[
\left(\wh{M}_1\times\wh{M}_2,\wh{\bsy{\theta}}_1\oplus\wh{\bsy{\theta}}_2\right),
\]
where $\wh{M}_i \cong M_i$ and locally $\wh{M}_1$ has coordinates
$$
x, q_1, q_1, q_2, q_3, q_4, a_1, a_2, a_3,
$$
while $\wh{M}_2$ has coordinates
$$
x, p_1, p_1, p_2, p_3, p_4, b_1, b_2, b_3.
$$
The Pfaffian systems $\wh{\bsy{\theta}}_1$ and $\wh{\bsy{\theta}}_2$ are locally generated as

$$
\begin{aligned}
\wh{\bsy{\theta}}_1= 
\Bigg\{da_1-a_2dx, dq_1-q_3dx, dq_2-\frac{a_1}{q_3}(q_1+q_2)dx, dq_3-&(a_1+a_3q_3)dx,\cr 
&dq_4-\frac{a_1q_4}{q_3}dx\Bigg\}\cr
\end{aligned}
$$
and
$$
\begin{aligned}
\wh{\bsy{\theta}}_2=
\Bigg\{db_1-b_2dy, dp_1-\frac{b_1}{p_4}(p_1+p_2)dy, dp_2-p_4dy, &dp_3-\frac{b_1p_3}{p_4}dy,\cr 
&dp_4-(b_1+b_3p_4)dy\Bigg\}.\cr
\end{aligned}
$$
As described in \cite{AFV}, the relationship between $\bsy{\theta},\ \wh{\bsy{\theta}}_1$ and $\wh{\bsy{\theta}}_2$ is that every integral manifold of $\bsy{\theta}$ is a superposition of integral manifolds of $\wh{\bsy{\theta}}_1$ and $\wh{\bsy{\theta}}_2$. The superposition formula is the map\footnote{See section 6 (appendix) for further details on the superposition formula; in particular, how it is defined and constructed.} 
\[
\bsy{\pi}:\wh{M}_1\times\wh{M}_2\to M
\] 
defined by
$$
\begin{aligned}
&\bsy{\pi}(x, \bsy{q}, \bsy{a}; y, \bsy{p}, \bsy{b})=\cr 
&\Bigg(\frac{-q_1+q_1p_3+q_4p_1-q_2+q_2p_3+q_1q_4}{q_4},\ \frac{q_4p_1+q_4p_2+q_2p_3-p_1-p_2+p_3p_2}{p_3},\cr 
&\ \ \ \ \ \ \ \ \ \ \ \ \ \ \ \ \ \ \ \ \ \ \ \ \ \ \      \frac{q_3(q_4-1+p_3)}{q_4},\ \frac{p_4(q_4-1+p_3)}{p_3},\ x,\ y,\ \bsy{a},\ \bsy{b}\Bigg)\cr
& \ \ \ \ \ \ \ \ \ \ \ \ \ \ \ \ \ \ \ \ \ \ \ \ \ \ \ \ \ \ \ \ \ \ \ \  \ \ \ \ \ \ \ \ \ \ \ \                                                                       =\Big(z_1,\ z_2, \ z_3,\ z_4, \ x,\ y,\ \bsy{a},\ \bsy{b}\Big) .\cr
\end{aligned}
$$
The usefulness of this factorisation of the integration problem for $\bsy{\theta}$ is not only that the integration of $\bsy{\theta}_i$ relies on ODE while that of $\bsy{\theta}$ relies on PDE, but also that the $\bsy{\theta}_i$ are locally equivalent to prolongations of the contact system on $J^1(\mathbb{R},\mathbb{R}^2)$. To see this, we turn to the characterisation of partial prolongations of such contact systems provided by \cite{Vassiliou2006a}, \cite{Vassiliou2006b}, which also provide simple procedures for finding the equivalence. To implement this, we compute the annihilators

\begin{eqnarray}
{\rm ann}\,\wh{\bsy{\theta}}_1=\wh{H}_1=\ \ \ \ \ \ \ \ \ \ \ \ \ \ \ \ \ \ \ \ \ \ \ \ \ \ \ \ \ \ \ \ \ \ \ \ \ \ \ \ \ \ \ \ \ \ \ \ \ \ \ \ \  \ \ \ \ \ \ \ \ \ \ \ \ \ \ \ \ \ \ \ \ \ \ \ \ \ \nonumber\cr
\cr
\left\{\P x+q_3\P {q_1}+\frac{a_1}{q_3}(q_1+q_2)\P {q_2}+(a_3q_3+a_1)\P {q_3}+\frac{q_4a_1}{q_3}\P {q_4}+a_2\P {a_1}, \P {a_2}, \P {a_3}\right\},\cr
\end{eqnarray}
\newpage
\begin{eqnarray}
{\rm ann}\,\wh{\bsy{\theta}}_2=\wh{H}_2=\ \ \ \ \ \ \ \ \ \ \ \ \ \ \ \ \ \ \ \ \ \ \ \ \ \ \ \ \ \ \ \ \ \ \ \ \ \ \ \ \ \ \ \ \ \ \ \ \ \ \ \ \  \ \ \ \ \ \ \ \ \ \ \ \ \ \ \ \ \ \ \ \ \ \ \ \ \ \nonumber\cr
\cr
\left\{\P y+\frac{b_1}{p_4}(p_1+p_2)\P {p_1}+p_4\P {p_2}+\frac{b_1p_3}{p_4}\P {p_3}+(b_3p_4+b_1)\P {p_4}+b_2\P {b_1}, \P {b_2}, \P {b_3}\right\}.\cr
\end{eqnarray}

The construction we have just described can be summarised by the following commutative diagram, where $\iota_i,\iota_2$ are inclusions and $\wh{\pi}_1, \wh{\pi}_2$ are projections such that $\wh{\pi}_j^*\bsy{\theta}_j=\wh{\bsy{\theta}}_j$.

\begin{center}
\begin{tikzpicture}
\draw [->, black](0,1.5) -- (0,-1.8); 
\draw (0,1.5) node[above]{$(\wh{M}_1\times\wh{M}_2,\wh{\bsy{\theta}}_1\oplus\wh{\bsy{\theta}}_2)$};
\draw [->,black](-0.1,1.45) -- (-1.5,0);
\draw (-2,0) node[below]{$(M_1,\bsy{\theta}_1)$}; 
\draw [->,black](0.1,1.45) -- (1.5,0);
\draw (2,0) node[below]{$(M_2,\bsy{\theta}_2)$};
\draw (0,-1.8) node[below]{$(M,\bsy{\theta})$};
\draw [->,black] (-1.6,-0.55) -- (-0.37,-1.8);
\draw [->,black] (1.6,-0.55) -- (0.37,-1.8);
\draw (0.25,-0.25)node[]{$\bsy{\pi}$};
\draw (1.2, -1.1)node[below]{${\iota_2}$};
\draw (-1.2, -1.1)node[below]{${\iota_1}$};
\draw (1.2, 1.1)node[below]{${\wh{\pi}_2}$};
\draw (-1.2, 1.1)node[below]{${\wh{\pi}_1}$};
\end{tikzpicture}
\end{center}

Locally, the maps $\wh{\pi}_1:\wh{M}_1\to M_1$ and $\wh{\pi}_2:\wh{M}_2\to M_2$ are given by
$$
\begin{aligned}
&\wh{\pi}_1(\bsy{p}, \bsy{a}, x)=(\bsy{z}, \bsy{a}, x),\ \ \wh{\pi}_2(\bsy{p}, \bsy{b}, y)=(\bsy{z}, \bsy{b}, y).
\end{aligned}
$$

We now show that each of the $\wh{H}_i$ is locally equivalent to 
the partial prolongation $C\langle 0,1,1\rangle$ of the contact distribution on $J^1(\mathbb{R},\mathbb{R}^2)$---that is, the contact distribution on $J^1(\mathbb{R},\mathbb{R}^2)$ partially prolonged so that one dependent variable has order 2 and the other order 3, with canonical local normal form 
\[
C\langle 0,1,1\rangle=\Big\{\P t+z^1_1\P {z^1}+z^1_2\P {z^1_1}+z^2_1\P {z^2}+z^2_2\P {z^2_1}+z^2_3\P {z^2_2},\ \ \P {z^1_2},\ \ \P {z^2_3}\Big\}.
\]
Let $\mathcal{D}$ be a smooth distribution on the manifold $M$, and assume that $\mathcal{D}$ is totally regular in the sense that $\mathcal{D}$, all its derived bundles, and all their corresponding Cauchy bundles have constant rank on $M$. Denote $m_i=\dim \mathcal{D}^{(i)}$, $\chi^i=\dim\, \ch\mathcal{D}^{(i)}$, and 
$\chi^i_{i-1}=\dim\,\ch\mathcal{D}^{(i)}_{i-1}$, where
\[
\ch\mathcal{D}^{(i)}_{i-1}=\mathcal{D}^{(i-1)}\cap\ch\mathcal{D}^{(i)}.
\] 
Below $k$ denotes the derived length of $\mathcal{D}$.

\vskip 5 pt
According to Theorem 4.1 of \cite{Vassiliou2006a}, a totally regular distribution $\mathcal{D}$ on smooth manifold $M$ is locally equivalent to a partial prolongation of the contact distribution on $J^1(\mathbb{R},\mathbb{R}^q)$ for some $q$ if and only if
\begin{enumerate}
\item The integers $m_i, \chi^j, \chi^j_{j-1}$ satisfy the numerical constraints

\begin{eqnarray}
\chi^j &=& 2m_j-m_{j+1}-1,\ \ \ 0\leq j\leq k-1\nonumber\cr
\chi^i_{i-1}&=& m_i-1,\ \ \ \ \ 1\leq i\leq k-1\cr
\end{eqnarray}

\item If $m_k-m_{k-1}>1$, then a certain canonically associated bundle called the {\it  resolvent} is integrable.\footnote{In the original formulation of Theorem 4.1 in \cite{Vassiliou2006a}, the integrability of $\ch\mathcal{D}^{(i)}_{i-1}$ is an additional hypothesis to be checked. This is a simple task, but unnecessary, since it is easy to see that this bundle is always integrable and that this hypothesis can be omitted.}
\end{enumerate}

A pair $(M,\mathcal{D})$ that satisfies these conditions is said to be a {\it Goursat manifold} or {\it Goursat bundle}. Moreover, if $\mathcal{D}$ is a Goursat bundle on $M$, then it is locally equivalent to a partial prolongation with $\chi^j-\chi^j_{j-1}$ dependent variables at order $j<k$ and $m_k-m_{k-1}$ dependent variables at highest order $k$. This uniquely identifies the partial prolongation associated to a given Goursat manifold. Before discussing the question of {\it constructing} equivalences, let us solve the recognition problem for the distributions $\wh{H}_i$. We will demonstrate the procedure for $\wh{H}_1$, which we denote temporarily by $\wh{K}$.

We find that
\[
\ch\wh{K}=\{0\},\ \ \ch\wh{K}^{(1)}=\{\P {a_2},\ \P {a_3}\},\ \ \ch\wh{K}^{(2)}_1=\{\P {a_2},\ \P {a_3},\ \P {a_1},\ \P {q_1}\} ,
\] 
\[
\ch\wh{K}^{(2)}=\{\P {a_2},\ \P {a_3},\ \P {a_1},\ \P {q_1},\ \P x\}.
\]
Calculation shows that the dimensions of the derived bundles are
\[
\dim\,\wh{K}=3,\ \dim\,\wh{K}^{(1)}=5,\ \dim\,\wh{K}^{(2)}=7,\ \dim\,\wh{K}^{(s)}=8,\ \ s\geq 3.
\]
Hence the derived length is $k=3$. Below we check the first condition of a Goursat bundle:
\vskip 10 pt
\begin{center}
\begin{tabular}{|c|c|c|c|c|c|}
\hline
$j$&$m_j$&$m_{j-1}-1$&$2m_j-m_{j+1}-1$&$\chi^j_{j-1}$&$\chi^j$\cr
\hline
0&3&$-$&$6-5-1=0$&$-$&0\cr
1&5&2&$10-7-1=2$&2&2\cr
2&7&4&$14-8-1=5$&$4$&5\cr
\hline
\end{tabular}
\end{center}
\begin{center}
Table: Checking the numerical constraints satisfied by $(M,\mathcal{D})$.
\end{center}

Therefore, $\wh{K}$ is a Goursat bundle with $k=3$, $m_k-m_{k-1}=1$, and the only nonzero difference $\chi^j-\chi^j_{j-1}$ being at order $j=2$: $\chi^2-\chi^2_1=5-4=1$. Hence there is one variable of order 2 and one variable of order 3. This solves the recognition problem, and we can assert that $\wh{K}:=\wh{H}_1$ is locally equivalent to $C\langle 0,1,1\rangle$. Next, we show how to construct an equivalence. Given a Goursat bundle, an efficient method for constructing an equivalence map was worked out in \cite{Vassiliou2006b}; it relies on the filtration induced on the cotangent bundle. Denote by $\Xi^{(j)}$ and $\Xi^{(j)}_{j-1}$ the annihilators of
$\ch\wh{K}^{(j)}$ and $\ch\wh{K}^{(j)}_{j-1}$, respectively. Then we obtain the filtration
$$
\Xi^{(2)}\subset\Xi^{(2)}_1\subset\Xi^{(1)},
$$ 
expressed explicitly as
$$
\{dq_1, dq_2, dq_4\}\subset\{dq_1, dq_2, dq_4, dx\}\subset\{dq_1, dq_2, dq_3, dx, dq_4, da_1\}.
$$
The construction proceeds by building appropriate differential operators and functions. Because  $m_k-m_{k-1}=1$, condition 2 of the definition of Goursat manifold is vacuous. Instead we fix any first integral of $\ch\mathcal{D}^{(k-1)}$, denoted $t$, and any section $Z$ of $\mathcal{D}$ such that $Zt=1$. Then, define distributions
$\Pi^k$ inductively as follows:
$$
\Pi^{\ell+1}=[Z,\Pi^\ell],\ \Pi^1=\ch\mathcal{D}^1_0,\ \ 1\leq \ell\leq k-1.
$$ 
There is a function $\varphi^k$ which is a first integral of $\Pi^k$ such that $d\varphi^k\wedge dt\neq 0$. The function $\varphi^k$ is said to be a {\it fundamental function of order} $k$.
The space of fundamental functions of lower order can be constructed from the filtration above by taking quotients. Specifically, as noted above, in this case the only fundamental functions of less than maximal order 3 are of order 2. They are described by the quotient bundle
$$
\Xi^{(2)}_1/\Xi^{(2)}=\{ dx \}  .
$$
Without loss of generality, we can take $\varphi^2=x$ to be a fundamental function of order 2. The construction of the equivalence map is now as follows. The function $t$ is the ``independent variable," and successive differentiation gives the higher order variables
$$
\z^1_0=\varphi^2,\ \z^1_1=Z\varphi^2,\ \z^1_2=Z^2\varphi^2,\ \z^2_0=\varphi^k,\ \z^2_1=Z\varphi^k,\ \z^2_2=Z^2\varphi^k,\ \z^2_3=Z^3\varphi^k. 
$$
We now implement this. The first integrals of $\ch\wh{K}^{(2)}$ are spanned by $q_1,q_2,q_3$, and any function of these can be chosen to be $t$. If we choose (say) $t=q_1$, then for $Z$ we choose
$$
Z=\frac{1}{q_3}X,
$$
where $X$ is the first vector field in the basis for $\wh{K}$ above; for then $Zt=1$, as required. We then construct the integrable distribution $\Pi^3$ as described above and discover that its first integrals are spanned by
$$
q_1,\ \ \frac{q_4}{q_1+q_2}.
$$
Hence the fundamental function of highest order (three) is
\[
\z^2_0=\varphi^k=\frac{q_4}{q_1+q_2}.
\]
The data
\[
t=q_1,\ \ \z^1_0=x,\ \ \z^2_0=\frac{q_4}{q_1+q_2}
\]
and differentiation by $Z$ now construct the local equivalence $\psi$ identifying $\wh{K}=\wh{H}_1$ and $C\langle 0, 1, 1\rangle$. The local inverse 
$\sigma_1=\psi^{-1}_1:\mathbb{R}\to \wh{M}_1$ determines the integral submanifolds of $\wh{H}_1$.

An exactly analogous calculation holds for $\wh{H}_2$, and one arrives thereby at an explicit map $\sigma_2=\psi^{-1}_2:\mathbb{R}\to\wh{M}_2$ representing the integral manifolds of $\wh{H}_2$. The explicit integral manifolds of $\bsy{\theta}$ are a {\it superposition} of those of $\wh{H}_1$ and $\wh{H}_2$:
\[
\mathbb{R}\times\mathbb{R}\to\bsy{\pi}(\psi^{-1}_1(\mathbb{R}),\psi^{-1}_2(\mathbb{R})).
\] 
The local diffeomorphisms $\psi_i$ are guaranteed to have a block-triangular structure and are therefore usually not difficult to invert in practice.

\vskip 15 pt
In this way we obtain remarkably compact representations for wave maps into this metric:

\vskip 5 pt

\begin{thm}[hyperbolic Weierstrass representation]\label{Wrep}
For each collection of twice continuously differentiable real valued functions  $k(s)$, $h(s)$, $m(t)$, $f(t)$ of parameters $s,t$, the functions
\begin{eqnarray}
x&=&f(t),\ \ \ \ \ \ \ \ y=h(s), \nonumber\cr\cr
u&=&-\frac{m\dot{k}s+mk+\dot{k}+k^2-tm\dot{k}}{m\dot{k}},\cr \cr
v&=&-\frac{m^2+tk\dot{m}+mk+\dot{m}-ks\dot{m}}{k\dot{m}}\cr
\end{eqnarray}
  
define harmonic maps
\[
(\mathbb{R}^{1,1},\ dx\,dy)\to \left(N, \frac{du_1^2-du_2^2}{2u_1}\right)\ \ \ \ \ \ {\rm by}\ \ \ \ \ \ u_1=u+v,\ \ \ u_2=u-v.
\]

\end{thm}

\section{Appendix: The superposition formula}
For completeness, in this appendix we make a remark on the construction of the superposition formula $\bsy{\pi} : \wh{M}_1\times\wh{M}_2\to M$.  A general construction valid for any decomposable exterior differential system was worked out in \cite{AFV}. Below we present results for the wave map system studied in this paper. The hyperbolic structure of the wave map system in adapted coordinates is given by
\[
H=H_1\oplus H_2,
\]
where
$$
H_1=\left\{\P x+z_3\P {z_1}+\frac{a_1}{z_3}(z_1+z_2)\P {z_2}+(a_3z_3+a_1)\P {z_3}+\frac{z_4a_1}{z_3}\P {z_4}+a_2\P {a_1}, \P {a_2}, \P {a_3}\right\},
$$

$$
H_2=\left\{\P y+\frac{b_1}{z_4}(z_1+z_2)\P {z_1}+z_4\P {z_2}+\frac{b_1z_3}{z_4}\P {z_3}+(b_3z_4+b_1)\P {z_4}+b_2\P {b_1}, \P {b_2}, \P {b_3}\right\}.
$$
Calculation shows that the infinitesimal symmetries of $H_1$ which are tangent to the level sets of all the first integrals $x, y, \bsy{a}, \bsy{b}$, known as the {\it tangential characteristic symmetries} of $H_1$,  are spanned by

$$
\mathfrak{e}_1=\left\{-z_4\P {z_2},-z_4\P {z_4},\frac{z_1+z_2}{z_4}\P {z_1}+\frac{z_3}{z_4}\P {z_3}+\P {z_4},\P {z_1}-\P {z_2}\right\}.
$$
Similarly, the tangential characteristic symmetries of $H_2$ are spanned by
$$
\mathfrak{e}_2=\left\{z_3\P {z_1}, z_3\P {z_3}, \frac{z_1+z_2}{z_3}\P {z_2}+\P {z_3}+\frac{z_4}{z_3}\P {z_4}, -\P {z_1}+\P {z_2}\right\}.
$$

The structure of these Lie algebras are in {\it  reciprocal relation}; namely, for $\mathfrak{e}_1$ the nonzero Lie brackets are
\[
[e_1,e_2]=-e_1,\ [e_1,e_3]=e_4,\ [e_2, e_3]=-e_3,
\]
while for $\mathfrak{e}_2$ we have
\[
[f_1,f_2]=f_1,\ [f_1,f_3]=-f_4,\ [f_2, f_3]=f_3,
\]
and
\[
[e_i,f_j]=0,\ \ \ \forall\ \ i,j.
\]

Recall the local product manifold  $\wh{M}_1\times\wh{M}_2$:
$$
\wh{M}_1\times\wh{M}_2=(\bsy{q},\bsy{a}, x;\bsy{p}, \bsy{b}, y)
$$
where
$$
\bsy{p}=p_1, p_2, p_3, p_4;\ \ \ \bsy{q}=q_1, q_2, q_3, q_4.
$$
The diagonal action $\mu_D:G\times \wh{M}_1\times \wh{M}_2\to \wh{M}_1\times \wh{M}_2$ takes the form
$$
\mu_D(g)(m_1,m_2)=\Big(\mu_1\big(g^{-1}\big)(m_1);\mu_2(g)(m_2)\Big),\ \ \ \forall\ \ g\in G,\ \ (m_1,m_2)\in\wh{M}_1\times\wh{M}_2,
$$
where $\mu_i$ is the action generated by $\wh{\pi}_i^*\mathfrak{e}_i$ on $\wh{M}_i$. This is the diagonal action mentioned in Theorem 2.2.

The Lie algebras $\mathfrak{e}_1,\mathfrak{e}_2$ define the isomorphism class of the {\it Vessiot algebra} for this Darboux integrable system, in the terminology of \cite{AFV}. Each $\mathfrak{e}_i$ is isomorphic to a solvable Lie algebra $\mathfrak{g}$, and each frames a neighbourhood of a point of $\mathbb{R}^4$. By the converse of Lie's second fundamental theorem, there is a local Lie group $(G, m)$, where $m: G\times G\to G$ denotes group composition, such that $\mathfrak{e}_i$ coincide with the infinitesimal left and right translations on $G$. From the expressions for either $\mathfrak{e}_1$ or $\mathfrak{e}_2$, by computing flows or otherwise, we can compute the function $m$. The components of $m$ coincide precisely with the first 4 components of the superposition map $\bsy{\s}$, where we interpret $\left(\{p_i\}_{i=1}^4,\{q_i\}_{i=1}^4\right)$ as local coordinates around $(e,e)\in G\times G$; $e$ being the identity in $G$. Thus, we see that the Cauchy characteristic vector field constructed in Theorem \ref{IVPmain} is a curve in the Lie algebra of tangential characteristic symmetries of $H_1$. There is an entirely analogous construction in which the first integrals $\alpha_j$ of $H_2$ are used instead of those of $H_1$, and in that case the Cauchy vector can be expressed as a curve in the tangential characteristic symmetries of $H_2$.






\bibliographystyle{model1b-num-names}



\end{document}